\documentclass{amsart}

\usepackage{graphicx}%
\usepackage{multirow}%
\usepackage{amsmath,amssymb,amsfonts}%
\usepackage{amsthm}%
\usepackage{mathrsfs}%
\usepackage[title]{appendix}%
\usepackage{xcolor}%
\usepackage{textcomp}%
\usepackage{manyfoot}%
\usepackage{booktabs}%
\usepackage{algorithm}%
\usepackage{algorithmicx}%
\usepackage{algpseudocode}%
\usepackage{listings}%
\usepackage{longtable}
\usepackage{nicematrix,tikz-cd}

\usepackage{natbib}
\setcitestyle{authoryear,open={(},close={)}} 
\newtheorem{theorem}{Theorem}[section]
\newtheorem{proposition}[theorem]{Proposition}
\newtheorem{lemma}[theorem]{Lemma}

\theoremstyle{definition}
\newtheorem{definition}{Definition}

\numberwithin{equation}{section}

\begin{document}
\title[Non-central distribution]{On non-central distribution of the matrix ratio}
\author[Haoming Wang]{Haoming Wang}
\address{Center for Combinatorics, Nankai University, Tientsin {\rm300071}, China}
\email{wanghm37@nankai.edu.cn}
\thanks{}
\subjclass[2020]{Primary {62H10}; Secondary {33C15}, {93B35}.}
\keywords{Matrix normal distribution, Non-central ratio distribution, Confluent hypergeometric function ${}_1F_1$ of a single matrix argument}
\date{}
\dedicatory{}

\begin{abstract}We derive the distribution of the ratio of a non-central mean matrix and a sample covariance matrix. This aligns with the confluent term ${}_1F_1$ in the non-central uni-variate Student's $t$. Some extensions of matrix-variate distributions are considered.
\end{abstract}

\maketitle

\section{Introduction}

In uni-variate statistical analysis, \cite{student1908} proved that if $x_1,x_2,\dots,x_N$ are i.i.d. sample values from $N(\mu,\sigma^2)$, then the sample covariance
\[s^2 = \sum (x - \bar{x})\]
is independent of the sample mean $\bar x = (1/n)\sum x$, by which one can define the Student $t$-distribution as the ratio of them, i.e.,
\[ t = \frac{\bar x}{s}.\]
This concerns with the uni-variate case. Albeit a bit difficult, the non-central distribution can be expressed through the confluent hypergeometric function ${}_1F_1$, which is sometimes hard to compute numerically.

Being aware of the complexity of this problem, \cite{Dunnett1954} generalised the uni-variate $t$-distribution to the bi-variate central case, which was then followed by \cite{John1961OnTE}, \cite{Kshirsagar1961SomeEO}, and \cite{Dickey1967} to the multivariate $t$-distribution afterwards, all three in the central case. As for the non-central scenario, existing literature on multi-variate and matrix-variate generalisation of $t$-distribution only considers the central ratio plus a fixed mean to avoid introducing the confluent term ${}_1F_1$. This can be seen from \cite{LIN1972339}, \cite{PHILLIPS1985157}, \cite{KotzNadarajah2004}, \cite{diaz2012matricvariate}, and \cite{gupta2018matrix}. 
However, in practice, the presence of location shifts or signal errors renders classical central approximations inaccurate, thereby distorting test power and confidence interval coverage. This motivates a more precise characterization of the non-central matrix ratio distribution and its density.

In this paper, we introduce the matrix $t$-distribution with a non-central matrix normal numerator in this ratio, in accordance with the non-central uni-variate $t$-distribution. 
However, if one insists on the most accurate assumption, the determination of the sample variances and covariances of $n$ samples in $p$ variables usually has, assuming normality for simplicity, $\frac{1}{2}np(np-1)$ correlation coefficients or regression coefficients to be estimated. A common approach to alleviate this burden is to assume the population covariance is governed by a specific tensor structure, that is, a suitable spectral decomposition like the elliptical assumption in \cite{1990Generalized}. In this sense, we only need to consider the distribution of $\frac{1}{2}p(p+1)$ quadratic forms,
\begin{equation*}
    \begin{aligned}
        {y}_1'{A}_{11}{y}_1, {y}_1'{A}_{12}{y}_2, \dots, {y}_1'{A}_{1p}{y}_p,\\
        {y}_2'{A}_{22}{y}_2, \dots, {y}_2'{A}_{2p}{y}_p, \\
        \vdots\qquad \\
        {y}_p'{A}_{pp}{y}_p,
    \end{aligned}
    \label{eq: illustration of triangular ararry heter}
\end{equation*}
which corresponds to four tensor forms $T_{1}$, $T_{1\frac{1}{2}}$, $T_{2}$, and $T_{3}$ of the population precision matrix that will be explained in the following.

\cite{khatri1966} first obtained the joint distribution of the homogeneous quadratic forms, that is, $A_{ij} = A$ a fixed covariance with the hypergeometric function ${}_0F_0$ of two matrix arguments. If furthermore, $A = I$, this reduces to the distribution of \cite{wishart1928generalised}. In Khatri's cited work, the exact distribution is not expressed in the series of zonal polynomials such as \cite{Kotz1967ab} but his method involves nuisance parameters $q_{ij} >0$. This expression is difficult in practice since the density is gamma-like in the central setting, where the hypergeometric function ${}_0F_0$ does not appear explicitly. In addition, hypergeometric functions of two matrix arguments often occurs in the moment generating function and distribution of latent roots for the central case, instead of the density as noted by \cite{james1960distribution, james1964distribution}, and they are usually computationally intractable. Thus, after nearly thirty years, the heterogeneous quadratic forms were reconsidered in \cite{Mathai1992} by calculating the cumulants only for $p=2$. It lies in another core subject of this article the determination of the central distribution of these $\frac{1}{2}p(p+1)$ heterogeneous quadratic forms for general $p>2$ using no hypergeometric function and unnecessary constant $q_{ij}$. To begin with a reasonable foundation, these matrices must satisfy $A_{ij}' = A_{ji}$ and the diagonal entries in $A_{ij}(i\neq j)$ are equal to zero. This promotes an ad hoc approach for the simultaneous diagonalisation that will be our main technique to handle this dimensionality curse.


This extended distribution, representing correlations in multi-variate normal samples, has been considered by \cite{dawid1977spherical}, \cite{kabe1984classical}, \cite{1990Generalized}, \cite{goodall1993multivariate}, \cite{CaroLopera2014OnGW}, \cite{gupta2018matrix}, \cite{mathai2022mul}. According to \cite{1990Generalized}, the spectral decomposition is hard to obtained if the samples are not independent and identically distributed. This paper utilizes an organised tensor forms $T_{1}$, $T_{1\frac{1}{2}}$, $T_{2}$, and $T_{3}$ of the population precision matrix to handle this problem. It directly assumes the spectral structure of the precision matrix and answer this question with the given stochastic representation. By this means, we derive the density of the matrix $t$-distribution with a similar confluent hypergeometric function ${}_1F_1$, which can be seen from the sensitivity analysis relatively close to the uni-variate Student's $t$-distribution.

\section{$T_{1}$, $T_{1\frac{1}{2}}$, $T_{2}$, and $T_{3}$}
To go with, let us introduce some basic knowledge of the spectral decomposition. Block an $np\times np$ Hermite matrix $\varTheta$ as
	\begin{equation*}
		\varTheta = \begin{bmatrix}
			\varTheta_{11} &\varTheta_{12} & \dots & \varTheta_{1n}\\
			\varTheta_{21} & \varTheta_{22} & \dots & \varTheta_{2n}\\
			\vdots & \vdots & \ddots & \vdots \\
			\varTheta_{n1} & \varTheta_{n2} & \dots & \varTheta_{nn}
		\end{bmatrix}
		\label{eq: block matrix}
	\end{equation*}
	where each $\varTheta_{ii}$ is positive definite of the order $p$. Assume $\varTheta$ is one of the following four types.
	\begin{itemize}
		\item There exist orthonormal column vectors $\{b_{j}\}$ of length $p$ such that $\varTheta$ can be developed into the sum
		\begin{equation*}
		    \varTheta = \sum_{j=1}^{n}\sum_{j^{\prime}=1}^{n} A_{jj^{\prime}} \otimes B_{jj^{\prime}}, \eqno{(T_{1})}
		\end{equation*}
		where $B_{jj^{\prime}} = b_{i} \overline{b_{j^{\prime}}}^{\prime}$ and $A_{jj^{\prime}}(k,l) = \overline{b_{j}}^{\prime} \varTheta_{kl} b_{j^{\prime}}$. In this case, $$A_{jj^{\prime}}(k,k)  = 0 \, (j \neq j^{\prime}).$$
        In a word, the diagonal entries of $A_{jj^{\prime}}$ vanish if $j\neq j'$.
		\item It is $T_{1}$ and additionally, 
        \begin{equation*}
            A_{jj^{\prime}} = A_{jj}\delta_{jj^{\prime}}, \eqno{(T_{1\frac{1}{2}})}
        \end{equation*}
        that is, $A_{jj^{\prime}}$ are vanishing if $j\neq j^{\prime}$.
		\item It is $T_{1}$ and there exist further orthonormal column vectors $\{a_{i}\}$ of length $n$ such that $\varTheta$ can be developed into the sum
		\begin{equation*}
		    \varTheta = \sum_{i=1}^{n}\sum_{j=1}^{p} \gamma_{ij} A_{i} \otimes B_{j}, \eqno{(T_{2})}
		\end{equation*}
		where $A_{i} = a_{i}\overline{a_{i}}^{\prime}$, $B_{j} = b_{j}\overline{b_{j}}^{\prime}$, and $\otimes$ is the Kronecker product operation.
		\item It is $T_{2}$ and additionally, there exist $\alpha_{i}, \beta_{j}$ such that $\gamma_{ij} = \alpha_{i} \beta_{j}$ or equivalently,
	\begin{equation*}
	    \varTheta =  \varPhi^{-1} \otimes \varPsi^{-1}, \eqno{(T_{3})}
	\end{equation*}	
        where $\varPhi^{-1}  = \sum_{i=1}^{n} \alpha_{i} a_{i} \overline{a_{i}}^{\prime}$ and $\varPsi^{-1}  = \sum_{j=1}^{p} \beta_{j} b_{j}\overline{b_{j}}^{\prime}$.
	\end{itemize}
Take simplicity for considerations, the following discusses real distributions.

\begin{definition}
    The $n \times p$ matrix population $X$ is said to follow the matrix normal distribution if $\operatorname{vec}(X')$, the vectorisation of $X'$, follows the multivariate normal distribution $N_{np}(0,\varSigma)$. 
\end{definition}

This proposition classifies the matrix normal populations according to four tensor decompositions of their precision matrices $\varTheta = \varSigma^{-1}$.

\begin{proposition}\label{prop: pdf} The probability density function $p(X)$ of a matrix normal population $X$ with its precision matrix $ \varTheta_{i}$ given by $T_{i}$, $i=1,1\frac{1}{2},2,3$ are

\begin{equation*}
                p(X) = \frac{| \varTheta_{1}|^{\frac{1}{2}}}{(2\pi)^{\frac{np}{2}}} 
            \operatorname{etr} \left(-\frac {1}{2} \sum_{j=1}^{p} { A}_{jj}  X{ B}_{jj}{ X}^{\prime} - \sum_{j=1}^{p}\sum_{j^{\prime}=j+1}^{p}{ A}_{jj^{\prime}}  X{ B}_{jj^{\prime}}{ X}^{\prime}\right),\eqno(T_1)
			\label{eq: T1 d}\\
\end{equation*}
\begin{equation*}
                p(X) = \frac{| \varTheta_{1\frac{1}{2}}|^{\frac{1}{2}}}{(2\pi)^{\frac{np}{2}}}\operatorname{etr} \left(-\frac {1}{2} \sum_{j=1}^{p} { A}_{jj} X{ B}_{jj}{ X}^{\prime}\right), \eqno(T_{1\frac{1}{2}})		\label{eq: T1.5 d}\\
\end{equation*}
\begin{equation*}
            p( X) =\frac{| \varTheta_{2}|^{\frac{1}{2}}}{(2\pi)^{\frac{np}{2}}}
             \operatorname{etr} \left(-\frac {1}{2} \sum_{i=1}^{n}\sum_{j=1}^{p} \gamma_{ij} { A}_{i} X{ B}_{j}{ X}^{\prime}\right),\eqno(T_2)
			\label{eq: T2 d}\\
\end{equation*}
\begin{equation*}
                p( X) = \frac{1}{(2\pi)^{\frac{np}{2}}| \varPhi|^{\frac{p}{2}}| \varPsi|^{\frac{n}{2}}}
               \operatorname{etr} \left(-\frac {1}{2}  \varPhi^{-1}  X  \varPsi^{-1} {{ X}}^{\prime}\right). \eqno(T_3)
				\label{eq: T3 d}
\end{equation*}
\end{proposition}

\begin{proposition}\label{prop: nested precision form 3}
$T_{1}\supset T_{1\frac{1}{2}} \supset T_{2} \supset T_{3}$. 
\end{proposition}

\begin{proposition} 
The sample covariance $ X^{\prime} X$, where $ X \sim T_{i}, i=1, 1\frac{1}{2}, 2, 3$ follows a matrix normal distribution, is positive definite with probability one if and only if $n  > p -1 $.
    \label{lem: dykstral70}
\end{proposition}

The proofs of these propositions can be found in \cite{wang2025}.

\section{Distribution of the sample covariance}

\begin{definition}
   The $p \times p$ matrix $(X+M)^{\prime} (X+M)$ is said to follow the product moment distribution if $ X \sim T_{i}, i=1, 1\frac{1}{2}, 2, 3$ is a matrix normal distribution, where $M$ is a fixed $n\times p$ matrix.
\end{definition}
In particular, if $X \in T_3$ and $\varPhi = I_n$, this definition reduces to the classical Wishart distribution.

\begin{theorem} \label{thm: main theorem}  Suppose $S \in T_1$ follows a product moment distribution. Assume that the diagonals in the $p\times p$ real matrix $ U = (u_{ij})$ with $u_{ij} = \operatorname{tr}( A_{ij})$ are all positive, that is, $u_{ii}>0 \, (i=1,2,\dots,p)$. Then the probability density distribution of $ S$ depends only on $T= (t_{ij})$ with $t_{ij} = \operatorname{tr}( B_{ij}  S)$ when $n > p - 1$,
        \begin{eqnarray}
         \frac{| \varTheta_{1}|^{\frac{1}{2}}}{2^{\frac{np}{2}}\Gamma_{p}(\frac{n}{2})} \operatorname{etr} \left(-\frac{1}{2}  U T \right) | T|^{\frac{n-p-1}{2}}          \text{ when } {M}=0; \\
                \times \operatorname{etr} \left(-\frac{1}{2}  \varOmega\right){}_{0}F_{1}\left(\frac{n}{2};\frac{1}{4} \varDelta  T\right) \text{ when } {M}\neq0;
        \end{eqnarray}
        where $\varOmega = \sum_{i,j=1}^{p} B_{ij}  M^{\prime}   A_{ij}  M$ and $ \varDelta = \sum_{i,j,k,l=1}^{p} B'_{ij} M'A'_{ij}A_{kl} M B_{kl}$.
\end{theorem}

To prove Theorem \ref{thm: main theorem}, we need three lemmas on the generalized hypergeometric function ${}_pF_q$.

\begin{lemma}[\cite{shimizu2022}]\label{lem: hypergeometric} Suppose $ A$ is an $n\times n$ and $ B$ is a $p\times p$ real symmetric matrix. If $n \geq p$,
    \[\begin{aligned}
        \int_{O(n)}{}_{0}F_{0}( A H_1 B H_1^{\prime}) [dH] = {}_{0}F_{0}( A, B),\\ 
         H = [ H_{1},  H_{2}], \text{ where }  H_1 \text{ is } n \times p,
    \end{aligned}\]
    where the integral runs over all $n\times n$ orthogonal matrices $O(n) = \{H:  H^{\prime} H =  I_{n}\}$, and $[dH]$ is the normalised Haar measure defined on $O(n)$.
\end{lemma}

\begin{lemma}[\cite{james1961zonal}]\label{lem: james55} Suppose $X$ is an $n\times p$ matrix.
    \[\begin{aligned}
        \int_{O(n)} \operatorname{etr}( X H_1^{\prime}) (d H) = {}_{0}F_{1}\left(\frac{1}{2}n;\frac{1}{4} X' X\right), \\ 
         H = [ H_{1},  H_{2}], \text{ where }  H_1 \text{ is } n \times p,
        \end{aligned}\]
    where the integral runs over the the orthogonal group $O(n)$.
\end{lemma}

	\begin{lemma}[\cite{khatri1966}] Let $ A$ be an $n \times n$ real symmetric matrix and $ B$ an $p \times p$ real symmetric positive definite matrix with $n \geq p$. Let $ X$ be an $n \times p$ real matrix. Then
		\begin{equation*}
			\int_{ X^{\prime}  X= S} \operatorname{etr} ( A  X  B  X^{\prime}) (d  X)= \frac{\pi^{\frac{np}{2}}}{\Gamma_{p}(\frac{n}{2} )} | S|^{\frac{n-p-1}{2}} {_0}F_0( A,  B S).\label{eq: Khatri1966}
		\end{equation*}\label{lem: Khatri66}
	\end{lemma}

Let $n,p$ be positive integers and $n\geq p$. Assume that $Z$ is a $n \times p$ real matrix with full column rank. Thus, $Z$ can be uniquely decomposed into $Z = HR^{\frac{1}{2}},$ where $H$ is an $n\times p$ matrix satisfying $H'H = I_p$ and $R$ is a $p\times p$ symmetric positive definite matrix. In particular, we can choose $H = Z(Z'Z)^{-\frac{1}{2}}$ and $R = Z'Z $.

The set of $n\times p$ matrices $H$ satisfying $H'H = I_p$ is said to be a Stiefel manifold, denoted by $V_{n,p}$. For any such $H\in V_{n,p}$, we can always expand it into an orthogonal matrix $K = (H,H_\perp) \in O(n)$, where the column vectors of $H_\perp \in V_{n,n-p}$ are orthogonal to the column vectors of $H\in V_{n,p}$. \cite{james1954} (Equation (8.19)) was the first to prove that
    \begin{equation*} \bigwedge_{i,j=1}^n dz_{ij} = 2^{-p} |Z'Z|^{\frac{1}{2}(n-p-1)} \bigwedge_{i\le j}^n d(z_i'z_j) \bigwedge_{i< j}^n h_i'dh_j \bigwedge_{i=k+1}^n \bigwedge_{j=1}^k h_i'dh_j, \label{eq: james1954}\end{equation*}
    where $H = (h_1,\dots,h_k), H_\perp = (h_{k+1},\dots,h_n)$. In matrix form,
\begin{lemma}\label{prop: herz55}
$(dZ) = 2^{-p} |R|^{\frac{1}{2}(n-p-1)} \cdot (dR) \cdot (dK)$, where $(dK) = {(H'dH)} \cdot {(H_\perp'dH)}$.
\end{lemma}
The geometric meaning of $(dK)$ is, \({(H'dH)} \) describes the rotation of the column vectors of the object \( H \) within the \( p \)-dimensional subspace it spans, thus corresponding to the volume element of $O(p)$, and \( {(H_\perp'dH)} \) describes the change of direction of the \( p \)-dimensional subspace in the \( n \)-dimensional space, that is, the movement of the subspace in the Grassmann manifold \( G_{n,p} \). The two together describe the Stiefel manifold \( V_{n,p} \) since the orthogonal group \( O(p) \) is the principal bundle of the fiber on the Grassmann $G_{n,p}$.
\begin{proof}[Proof of Theorem \ref{thm: main theorem}]
    Suppose the central part holds for $n > p-1$. Decomposing $ X =  H  Z$ where $ Z$ is upper-triangular and $ H' H =  I_p$.
    Let $K = (H,H_{\perp}) \in O(n)$. From Lemma \ref{prop: herz55}, we could rewrite \eqref{eq: T1 d} as 
    \[\begin{aligned}
        \operatorname{etr} \left(-\frac{1}{2}  \varOmega\right)\int_{O(n)}\operatorname{etr}\left(-\frac{1}{2} \sum_{i,j=1}^{p}( A_{ij} M B_{ij})Z'H' \right) (d K) \\
        = \operatorname{etr} \left(-\frac{1}{2}  \varOmega\right){}_{0}F_{1}\left(\frac{n}{2};\frac{1}{4}\sum_{i,j,k,l=1}^{p}( B'_{ij} M'A'_{ij}A_{kl} M B_{kl}) Z'Z\right).    \end{aligned}\]
    Thus, if the central part is true, Theorem \ref{thm: main theorem} is consequently proved by applying Lemma \ref{lem: james55}. 
    
    From Lemma \ref{lem: Khatri66} on the quadratic forms in normal vectors, we can reduce the central part of Theorem \ref{thm: main theorem} to the $\frac{1}{2}p(p+1)$ independent quadratic forms by introducing $ Y =  X B = ( y_{1}, y_{2}\dots,  y_{p})$
    \begin{equation*}
    \begin{aligned}
        {y}_1'{A}_{11}{y}_1, {y}_1'{A}_{12}{y}_2, \dots, {y}_1'{A}_{1p}{y}_p,\\
        {y}_2'{A}_{22}{y}_2, \dots, {y}_2'{A}_{2p}{y}_p, \\
        \vdots\qquad \\
        {y}_p'{A}_{pp}{y}_p,
    \end{aligned}
\end{equation*}
In fact, we have for each ${y}_i'{A}_{ij}{y}_j$ the contribution to the probability density function
\[\begin{aligned}
    {_{0}}F_{0}\left( I - \frac{q_{ij}}{2} A_{ij}, q_{ij}^{-1} y_{i} y_{j}^{\prime}\right) = {_{0}}F_{0}\left( I - \frac{q_{ij}}{2} A_{ij}, q_{ij}^{-1} y_{j}^{\prime} y_{i}\right).\end{aligned}\]
However, we may find $ y_{j}^{\prime} y_{i} = t_{ij} = \operatorname{tr} ( B_{ij} S)$ so that from these properties of hypergeometric functions
\[\begin{aligned}
    {}_{0}F_{0} ( X,c Y) = {}_{0}F_{0} (c X, Y)\\
    {}_{0}F_{0} ( X,  I) = {}_{0}F_{0} ( X) = \operatorname{etr}( X), 
\end{aligned}\]
the terms concerning $q_{ij}$ cancel out. This yields the desired form.
\end{proof}

This extends the separable covariance $ \varPhi \otimes  \varPsi$ with potentially variable-level ($ \varPsi$) and/or sample-level ($ \varPhi$) correlations in the classical case, especially the result of \cite{khatri1966}. For example, assume the central part of Theorem \ref{thm: main theorem} holds and $ A_{ij} = \beta_{j}\delta_{ij}   \varPhi^{-1}$ such that $ \varPsi^{-1} = \sum_{j=1}^{p}\beta_j b_{j} b_{j}^{\prime}$, where $\beta_{j},  b_{j}$ are the latent roots and vectors of $ \varPsi^{-1}$. This implies that 
\[\begin{aligned}
    \operatorname{etr}\left(-\frac{1}{2} U  T \right) = \operatorname{etr}\left(- \frac{1}{2}\operatorname{tr} ( \varPhi^{-1})\sum_{j=1}^{p}\beta_{j}t_{jj}\right)     = \operatorname{etr}\left(- \frac{1}{2} \varPhi^{-1}\operatorname{tr} ( \varPsi^{-1}  S)\right) \\
    = \operatorname{etr} \left(- q^{-1} \varPsi^{-1}  S\right){}_{0}F_{0} \left( I - \frac{q}{2} \varPhi^{-1}, q^{-1} \varPsi^{-1}  S\right),
\end{aligned}\]
where, in the second equality, we used the definition $ T =  B  S  B'$. Thus,  Lemma \ref{lem: Khatri66} is proved.

\section{Student's $t$-distribution}

\begin{definition}
    The $m\times p$ matrix $Z = (X+M) S^{-\frac{1}{2}}$ is said to follow the matrix $t$-distribution if $X \in N_{m,p}(0,I_m,I_p)$ follows a matrix normal distribution and $S = Y'Y$ follows an independent product moment distribution, where $Y$ is a central $T_1$ distribution and $M$ is an $m\times p$ fixed matrix.
\end{definition}

\begin{theorem}\label{thm: matrix t} The central matrix $t$-distribution has the density
\begin{equation}\begin{aligned}
    \frac{\Gamma_p\left(n + \frac{1}{2}m -\frac{1}{2}\right)}{(2)^{\frac{1}{2}(n+m)p}\pi^{\frac{1}{2}mp}\Gamma_p(\frac{n}{2})}|U+Z'Z|^{- n - \frac{m}{2} -\frac{1}{2}},
\end{aligned}
\end{equation}
and the noncentral distribution is the central one multiplied by
\begin{equation}\begin{aligned}
\operatorname{etr} \left(- \frac{1}{2}M'M \right){}_1F_1\left( n+\frac{1}{2}m-\frac{1}{2},\frac{1}{2}p;\frac{1}{4}M'ZZ'MB(U+Z'Z)^{-1}B^{'}\right).
\end{aligned}
\end{equation}
\end{theorem}


        \begin{lemma}\label{lem1} Let $C$ be a $p\times p$ real symmetric positive definite matrix, and $D$ be a $p \times p$ arbitarty fixed matrix.
        \begin{equation}\begin{aligned}
           \int_{X>0} \operatorname{etr}& (-CX'X + DX')|X'X|^{a - \frac{p+1}{2}} (d X) \\
                    & = {\Gamma_p\Big(a-\frac{1}{2}\Big)}{|C|^{-(a-\frac{1}{2})}}{}_1F_1\left(a - \frac{1}{2};\frac{1}{2}p;\frac{1}{4}D'DC^{-1} \right)
        \end{aligned}
        \label{eq: gamma integral}
        \end{equation}
        where $\Re(a) >\frac{1}{2}(p-1)$, and $X$ runs over all $p\times p$ p.d. matrices.\end{lemma}
        
    \begin{proof}[Proof of Lemma \ref{lem1}] 
    Put $Y = XC^{\frac{1}{2}}$ and $E = DC^{-\frac{1}{2}}$. By Lemma \ref{prop: herz55}, any $n\times p$ matrix $Y$ uniquely decomposes as $Y = HR^{\frac{1}{2}}$, where $H = Y(Y'Y)^{-\frac{1}{2}}, R = Y'Y$ such that the $n \times p$ matrix $H$ belongs to the Stielfel manifold $V_{n,p}$. Let $K = (H,H_{\perp}) \in O(n)$. Thus, \eqref{eq: gamma integral} becomes 
    \[2^{-p}|C|^{-a+\frac{1}{2}}\int_{R>0}(dR)\int_{O(n)} \operatorname{etr}(- R + H'ER^{\frac{1}{2}})|R|^{a - \frac{p}{2}-1} (dK).\]
    By Lemma \ref{lem: james55}, this integral reduces to 
    \[2^{-p}|C|^{-a+\frac{1}{2}}\int_{R>0} |R|^{a - \frac{p}{2}-1} \operatorname{etr}(-R) {}_0F_{1} \left(\frac{1}{2}n;\frac{1}{4} E'ER\right) (d R).\]
    From \cite{constantine1963some}'s lemma,    
    we could expand ${}_0F_{1} \left(\frac{1}{2}n; E'ER\right) $ in terms of zonal polynomials to derive the confluent hypergeometric function ${}_1F_{1}$ term.
    \end{proof}

\begin{proof}[Proof of Theorem \ref{thm: matrix t}]
    The joint density of $S$ and $X$ is given by 
    \[\frac{1}{(2)^{\frac{1}{2}(n+m)p}\pi^{\frac{1}{2}mp}\Gamma_p(\frac{n}{2})} \operatorname{etr} \left(- \frac{1}{2}X'X-\frac{1}{2}  U T \right) | T|^{\frac{n-p-1}{2}}.\]
    Now, let $Z = (X + M)S^{-\frac{1}{2}}$. The Jacobian of the transformation is $|S|^{\frac{m}{2}}$. Substituting for $X$ in terms of $Z$ in the joint density of $X$ and $S$, and multiplying the resulting expression by $|S|^{\frac{m}{2}}$, we get the joint density  of $Z$ and $S$ as
    \[\begin{aligned}
        \frac{\pi^{-\frac{1}{2}mp}}{(2)^{\frac{1}{2}(n+m)p}\Gamma_p(\frac{n}{2})} \operatorname{etr} \left(- \frac{1}{2}M'M + S^{\frac{1}{2}}Z'M-\frac{1}{2}  (U + Z'Z)B'SB \right)
        |S|^{\frac{n+m-p-1}{2}},
    \end{aligned}\]
    where we have substituted $S = B'TB$ by definition. 
    By Lemma \ref{lem1}, we can derive the density of $Z$ by integrating $S^{\frac{1}{2}}$.
\end{proof}

\section{Sensitivity analysis}
Consider these two models
\[\text{(A)} \ Z = (X+M)S^{-\frac{1}{2}}, \qquad \text{(B)} \ Z =XS^{-\frac{1}{2}}+M.\]
We are going to analyse the sensitivity using 1st order moments
\[\frac{\partial E[Z_A]}{\partial M} = \frac{\partial E[MS^{-{1}/{2}}]}{\partial M} = E[S^{-1/2}],\]
\[\frac{\partial E[Z_B]}{\partial M} = \frac{\partial E[M]}{\partial M} = I.\]
If $S \sim W_{n,p}(0,\Sigma)$, then 
\[E[S^{-1}] = \frac{\Sigma^{-1}}{n-p-1}.\]
There is no simple closed-form general expression of $E[S^{-1/2}]$ (the square root of a matrix and its expectation are not easily commutative). However, it can be obtained numerically or using Monte-Carlo methods. Define
\[\delta = \|E[S^{-1/2}] - I\|_{F}, \quad r_{A} = \frac{\|E[MS^{-1/2}]\|_{F}}{\|M\|_F},\]
where $\|\cdot\|_F$ is the Frobenius norm.

\begin{figure}[!ht]	
	\begin{minipage}[h]{.49\linewidth}
            \centering
            \includegraphics[width=1\linewidth]{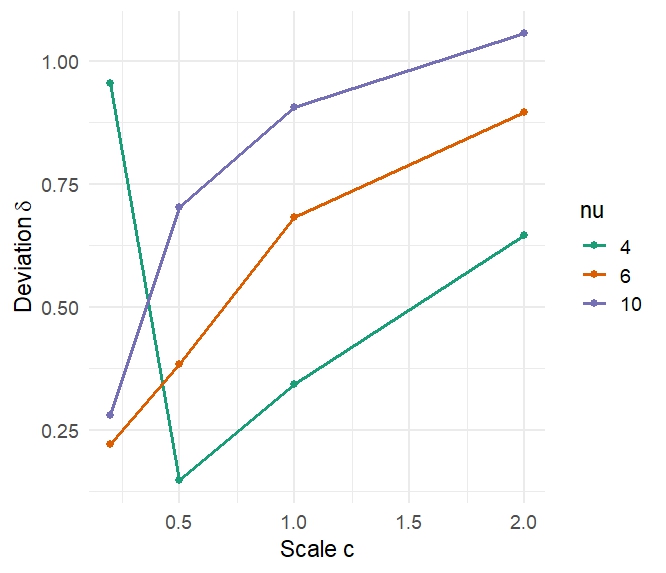}
    \end{minipage} 
	\hfill    
	\begin{minipage}[h]{.49\linewidth}
        \centering
        \includegraphics[width=1\linewidth]{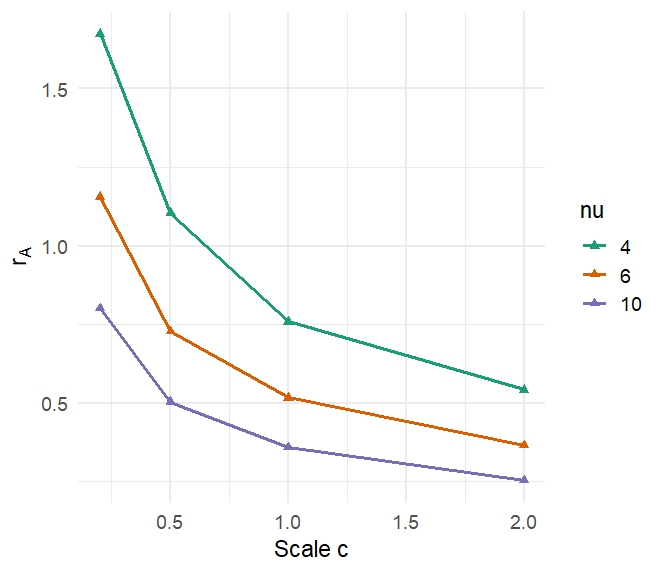}
	\end{minipage}
            \caption{Deviation $\delta$ and $r_A$ vary with different scales $c$ in $\Sigma = c I_{2}$ and $n=4,6,10$ with a given mean $M=[0,1;0,0]$.}
\end{figure}

\begin{figure}[!ht]	
	\begin{minipage}[h]{.49\linewidth}
		\centering
		\includegraphics[width=1\linewidth]{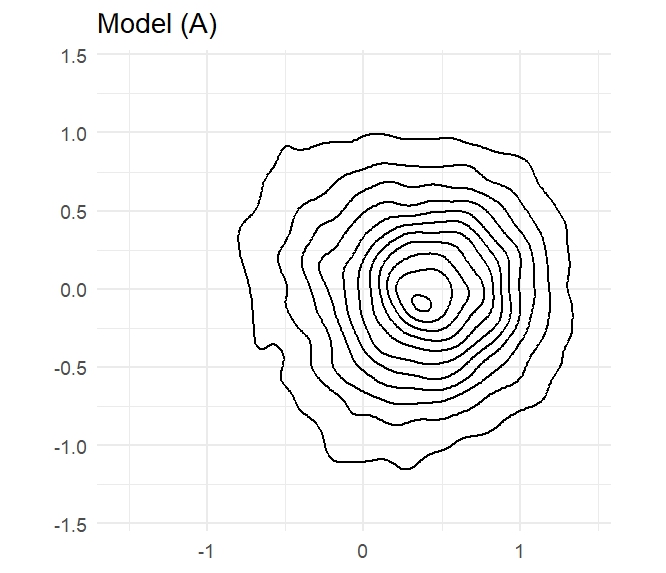}
	\end{minipage} 
	\hfill    
	\begin{minipage}[h]{.49\linewidth}
		\centering
		\includegraphics[width=1\linewidth]{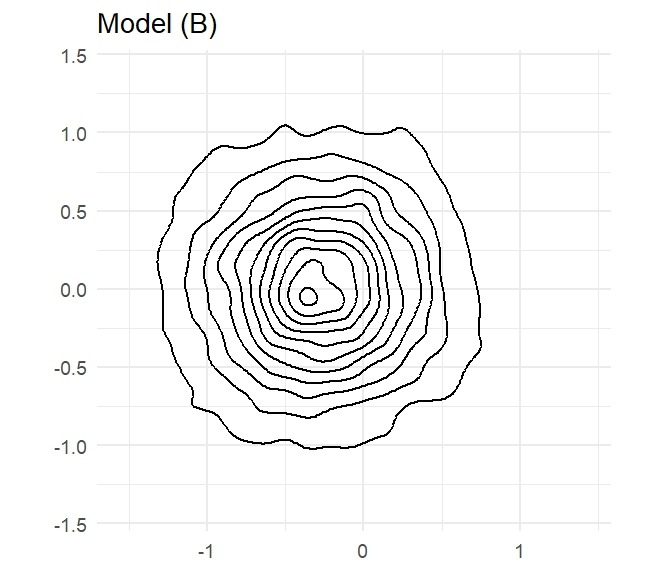}
	\end{minipage}
	\caption{2-D Counter of Student's $t$. Model (A) on the left and Model (B) on the right, where $X \sim N_{2,2}(0,I_2,I_2)$, $S\sim W_{2,2}(0,I_2)$, and $M = (1.0, 0.2; -0.3, 0.5)$.}
\end{figure}

\bibliographystyle{spbasic}
\bibliography{Sample}
\end{document}